\documentclass[letterpaper, 10 pt, conference]{IEEEtran}
\IEEEoverridecommandlockouts
\usepackage{cite}
\usepackage{amsmath,amssymb,amsfonts}
\usepackage{algorithmic}
\usepackage{import}
\usepackage{graphicx}
\usepackage{textcomp}
\usepackage{xcolor}
\usepackage{hyperref}
\usepackage{mathtools}
\usepackage{enumitem}

\usepackage{lipsum}

\usepackage{tikz}
\usetikzlibrary{positioning, calc,arrows.meta}

\usepackage[absolute,overlay]{textpos} 

\usepackage{cuted}
\usepackage{dblfloatfix} 

\newtheorem{thm}{Theorem}

\newtheorem{rem}{Remark}
\newtheorem{lem}{Lemma}

\newcommand{\Rank}[1]{\mathrm{rank}\left(#1\right)}

\newcommand{\Cauchy}[1]{\mathcal{C}(#1)}
\newcommand{\D}{\mathrm{d}}
\newcommand{\Lie}{\mathrm{L}}

\newcommand{\X}{\mathcal{X}}

\newcommand{\TX}{\mathcal{T}(\mathcal{X})}

\newcommand{\Span}[1]{\mathrm{span}\left\{#1\right\}}
\newcommand{\crksub}[1]{\underset{#1}{\subset}}

\newcommand{\xdot}{\dot{x}}

\newcommand{\zdot}{\dot{z}}

\newcommand{\zb}{\bar{z}}

\newcommand{\pad}[1]{\partial_{#1}}

\newcommand{\Ld}{\mathcal{L}}

\newcommand{\xidot}{\dot{\xi}}
\newcommand{\chidot}{\dot{\chi}}
\newcommand{\zetadot}{\dot{\zeta}}

\newcommand{\Dist}{\mathcal{D}}
\newcommand{\E}{\mathcal{E}}
\newcommand{\F}{\mathcal{F}}

\newcommand{\hook}{\mathbin{\lrcorner}}

\newcommand{\rotsubset}{\mathrel{\rotatebox[origin=c]{90}{$\underset{1}{\subset}$}}}

\newcounter{MYtempeqncnt}

\begin{document}

\title{ A Flat Triangular Structure Based on  \\a Multi-Chained Form \thanks{This research was funded in whole, or in part, by the Austrian Science Fund (FWF) P36473. For the purpose of open access, the author has applied a CC BY public copyright licence to any Author Accepted Manuscript version arising from this submission.}
}

\author{Georg Hartl, Conrad Gstöttner and Markus Schöberl \thanks{All authors are with the Institute of Control Systems, Johannes Kepler University Linz, Altenberger Strasse 69, 4040 Linz, Austria, \mbox{E-mail:}~{\tt \{georg.hartl, conrad.gstoettner, markus.schoeberl\}@jku.at}.
}}

\maketitle

\begin{abstract}
Determining whether a nonlinear multi-input system is differentially flat remains challenging. One way to obtain computationally tractable sufficient conditions is to give complete characterizations of flat normal forms. We introduce a structurally flat triangular form for control-affine systems with at least three inputs that is based on a multi-chained form. For two specific instances of this structure, we provide complete geometric characterizations, i.e., necessary and sufficient conditions under which a control-affine system is static-feedback equivalent to the respective triangular form. These characterizations yield sufficient conditions for differential flatness and, in turn, constructive procedures for computing flat outputs.
\end{abstract}


\section{Introduction}
Within nonlinear control theory, a central objective is to identify coordinate-free system properties that allow an equivalent representation in a canonical form---most notably, the integrator chain---under state and input transformations. Among the concepts that emerged from this pursuit are static and dynamic feedback linearization, leading to the notion of differential flatness, formulated in the early 1990s \cite{fliess_flatness_1995}. 

Formally, a control system of the form
\begin{equation}\label{eq:gen_nonlin_sys}
    \xdot = f(x,u)
\end{equation}
with $n$ states $x$ and $m+1$ control inputs $u$, is said to be \emph{flat} if there exist $m+1$ differentially independent functions
\begin{equation}\label{eq:intr_flat_ouput}
    \varphi(x, u, u^{(1)}, \ldots, u^{(\nu)}),
\end{equation}
where $u^{(\nu)}$ denotes the $\nu$-th time derivative of $u$, such that
\begin{equation}\label{eq:intr_flat_para}
    \begin{aligned}
        x & = F_x(\varphi, \varphi^{(1)}, \ldots, \varphi^{(r-1)}), \\
        u & = F_u(\varphi, \varphi^{(1)}, \ldots, \varphi^{(r)}) \; 
    \end{aligned}
\end{equation}
holds identically. In that case, \eqref{eq:intr_flat_ouput} is a \emph{flat output} of \eqref{eq:gen_nonlin_sys} and \eqref{eq:intr_flat_para} is its associated \emph{flat parameterization}. Consequently, every state and input trajectory can be parameterized by a trajectory \(\varphi(t)\) and finitely many of its time derivatives, thereby enabling systematic methods for trajectory planning and tracking \cite{fliess_lie-backlund_1999}. Within a differential–geometric framework, recent work has developed systematic procedures for designing tracking controllers for \((x,u)\)-flat systems, i.e., systems that admit a flat output \(\varphi(x,u)\) \cite{gstottner_tracking_2024}. If \eqref{eq:gen_nonlin_sys} admits a flat output depending only on the state, we say the system is \emph{\(x\)-flat} and \(\varphi(x)\) is an \emph{\(x\)-flat output}.

Given the fact, that no computationally tractable necessary and sufficient conditions for flatness of general nonlinear multi-input systems have yet been established, determining a flat output remains a non-trivial task. However, by focusing on subclasses of differentially flat systems, tremendous progress has since been made. For recent works see, e.g., \cite{schoberl_implicit_2014, nicolau_flatness_2017, gstottner_finite_2021, gstottner_necessary_2023,levine_differential_2025, nicolau_dynamic_2025}. 

We consider $n$th-order control-affine systems of the form
\begin{equation}\label{eq:ai_sys_m}
    \dot{x} = f(x) + g_0(x)u^0 + \cdots + g_m(x)u^m,
\end{equation}
where \(f\) is the \emph{drift vector field} and \(g_i\) are the \emph{input vector fields}. When the drift $f$ is absent, we call \eqref{eq:ai_sys_m} a driftless system.

Among approaches to the analysis of differentially flat systems, structurally flat triangular forms provide an effective framework for determining flat outputs, which typically appear as top variables in such normal forms. A driftless control-affine system with two inputs is flat if and only if it is (locally) static feedback equivalent to a system of the form
\begin{equation*}
    \zdot^1=w^0,\zdot^2=z^3w^0,\ldots,\zdot^{n-1}=z^nw^0,\zdot^n=w^1,
\end{equation*}
which is called the \textit{chained form} \cite{martin_feedback_1994}, see also \cite{murray_nilpotent_1994}. In  \cite{li_characterization_2013}, this normal form was extended by a drift compatible with the triangular structure yielding the \textit{extended chained form}. More general triangular forms presented in \cite{bououden_triangular_2011, gstottner_flat_2021, gstottner_structurally_2022, gstottner_triangular_2024, hartl_triangular_2025} were motivated by the extended chained form and enable constructive procedures for the computation of flat outputs.

Our focus, however, is on control-affine systems \eqref{eq:ai_sys_m} with at least three inputs (\(m+1\ge 3\)).
In the driftless case, a structurally flat triangular form is given by the \mbox{\emph{extended Goursat normal form}} \vspace{-1em}
\begin{equation*}
\arraycolsep=1pt
\begin{array}{rclrclcrcl}
    \zdot^0&=&w^0,&\zdot^1_1 & = & z^2_1w^0, &  & \zdot^1_m & = & z^2_mw^0, \\
    &&&& \vdots & & \cdots & & \vdots & \\
    &&&\zdot^{k_1-1}_{1} & = & z^{k_1}_{1}w^0, &  & \zdot^{k_m-1}_{m} & = & z^{k_m}_{m}w^0, \\
        &&& \zdot^{k_1}_1 & = & w^1, &  &  \zdot^{k_m}_m & = & w^m, \\
\end{array}
\end{equation*}
which is a special case of a multi-chained form \cite{murray_nonholonomic_1993, murray_nilpotent_1994}. Restricting further to chains of equal length (i.e., \mbox{$k_1=\cdots=k_m=k$}) yields the \emph{canonical contact form} that admits a complete geometric characterization \cite{respondek_transforming_2001, RespondekCanonicalContactSystems2002}. 
Allowing a drift that preserves this triangular structure leads to the \emph{canonical contact form with compatible drift} \cite{LiMultiinputControlaffineSystems2016}
\begin{equation}\label{eq:CCFD} 
    \arraycolsep=1.4pt
    \begin{array}{rclcrcl}
        \zdot^0 & = & w^0, &&&& \\
        \zdot^1_1 & = & z^2_1w^0 + a^1_1(z^0, \bar{z}_2), & & \zdot^1_m & = & z^2_mw^0 + a^1_m(z^0, \bar{z}_2), \\
        & \vdots & & \cdots & & \vdots & \\
        \zdot^{k-1}_1 & = & z ^k_1w^0 + a^{k-1}_1(z), & & \zdot^{k-1}_m & = & z ^k_mw^0 + a^{k-1}_m(z), \\
        \zdot ^k_1 & = & w^1, & & \zdot ^k_m & = & w^m, \;
    \end{array}
\end{equation}
with \mbox{$\zb_j = (z^1_1, \ldots, z^1_m, \ldots, z^j_1, \ldots, z^j_m)$} for \mbox{$1\leq j \leq k$}.  

Inspired by the works \cite{gstottner_flat_2021, gstottner_structurally_2022} about systems with \mbox{$m+1=2$} inputs, where the extended chained form is augmented by input prolongations and integrator chains appended to the top states, we extend these approaches to systems with $m+1\ge 3$ control inputs. 
First, we place integrator chains of equal length before the $m$ chained inputs and add a chain shorter by $s$ before the distinguished input $w^0$. Second, we append integrator chains of arbitrary length to the top state variables of \eqref{eq:CCFD}. Subsequently, we present necessary and sufficient geometric conditions for a system \eqref{eq:ai_sys_m} to be static feedback equivalent to two special cases of the presented triangular form, i.e., for $s=0$ and for $s=1$. As an application, we show that the gantry crane in three dimensions is static feedback equivalent to the case $s=1$, and, leveraging the characterization, we provide a constructive procedure to compute a flat output of this system. Our work is organized as follows: We introduce notation and terminology in Section \ref{sec:notation} and summarize relevant known results in Section \ref{sec:known_results}. Our main contributions are then presented in Section \ref{sec:main_results}. Finally, Section \ref{sec:example} illustrates our main theorems on a practical example. Sketches of the proofs of the main theorems are given in the appendix.
\vspace{-1.5ex}
\section{Notation and Terminology}
\label{sec:notation}

We adopt tensor notation and the Einstein summation convention. We omit index ranges when clear from the context. Let \(\X\) be an \(n\)-dimensional smooth manifold with local coordinates \(x=(x^1,\ldots,x^n)\) and tangent bundle \(\TX\). For an \(m\)-tuple of smooth functions \mbox{\(h=(h^1,\ldots,h^m):\X\to\mathbb{R}^m\)}, its Jacobian is denoted \(\partial_x h\in\mathbb{R}^{m\times n}\) and we denote \(\frac{\partial h^j}{\partial x^i}\) shortly by $\partial_{x^i}h^j$. The differentials \((\D h^1,\ldots,\D h^m)\) are abbreviated by \(\D h\) and the \(k\)-fold Lie derivative of a scalar function \(h^j\) along a vector field \(v\) is denoted by \(\Lie_v^k h^j\).

Let $v,w$ be two vector fields, $[v,w]$ indicates their Lie bracket. Further, let \(\Dist_1,\Dist_2\) be two distributions. Then, \([\Dist_1,\Dist_2]\) is the span of all brackets \([v_1,v_2]\) with \(v_i\in \Dist_i\). Similarly, \([v,\Dist_1]\) is the span of all brackets \([v,w_j]\) with \(w_j\in \Dist_1\). Given a distribution \(\Dist\), its derived flag is
\[
  \Dist^{(0)}=\Dist, \quad \Dist^{(i+1)}=\Dist^{(i)}+[\Dist^{(i)},\Dist^{(i)}] \;\; \text{for} \;\; i\ge 0. 
\]
The involutive closure \(\bar \Dist\) is the smallest involutive distribution containing \(\Dist\). The Cauchy characteristic distribution of $\Dist$ is defined by $\Cauchy{\Dist}=\{\,v\in \Dist\mid [v,\Dist]\subset \Dist\,\}$. That means, it is spanned by all vector fields $v \in \Dist$ for which the Lie bracket with any vector field in $\Dist$ remains in $\Dist$. Throughout, \(\subset\) is used inclusively allowing equality. The notation \(\Dist_1\underset{k}{\subset}\Dist_2\) indicates that the corank of \(\Dist_1\) in \(\Dist_2\) is \(k\).

Consider a control-affine system of the form \eqref{eq:ai_sys_m}. Geometrically, such a system is represented by the drift vector field $ f = f^i(x)\pad{x^i}$ and the input vector fields $g_j = g^i_j\pad{x^i}$ on the state manifold $\X$. We call two control-affine systems \emph{static feedback equivalent} (SFE), if they are locally equivalent via a diffeomorphism $z=\Phi(x)$ with its inverse $x=\hat{\Phi}(z)$ and a static input transformation $w^j = \alpha^j(x) + \beta^j_k(x)u^k$ with its inverse \mbox{$u^k = -\hat{\beta}^k_j(x)\alpha^j(x) + \hat{\beta}^k_j(x)w^j$}. The equivalent system reads\vspace{-1ex}
\begin{equation*}
    \zdot = a(z) + b_0(z)w^0 + \cdots + b_m(z)w^m \; ,
\end{equation*}
where $a^i(z)=(\pad{x^l}\Phi^i(x)(f ^l(x) - g ^l_k(x)\hat{\beta}^k_j(x)\alpha^j(x)))\circ\hat{\Phi}(z)$ and $b^i_j(z) = (\pad{x^l}\Phi^i(x)g ^l_k(x)\hat{\beta}^k_j(x))\circ\hat{\Phi}(z)$. 

Throughout, all functions, vector fields, and covector fields are assumed to be smooth, and all distributions and codistributions have locally constant rank. Unless stated otherwise, we consider generic points only.

\begin{figure*}[!b]
\normalsize
\setcounter{MYtempeqncnt}{\value{equation}}
\setcounter{equation}{5}
\vspace*{-4pt}
\begin{subequations}\label{eq:TFCCS+IC_q}
    \begin{equation}\label{eq:TFCCS+IC_xi}
    \hspace{-9.8em}\Sigma_\xi:\left\{ \hspace{1.1em}
      \begin{array}{rclcrclcrcl}
        \dot{\xi}^1_{0} &=& \xi^2_{0} & &
        \dot{\xi}^1_{1} &=& \xi^2_{1} & &
        \dot{\xi}^1_{m} &=& \xi^2_{m} \\
        & \vdots & & \hspace{1.15em} & & \vdots & \hspace{5.62em} & \cdots & \hspace{5.62em} & \vdots & \\
        \dot{\xi}^{k^{\xi}_{0}}_{0} &=& \chi^0 & &
        \dot{\xi}^{k^{\xi}_{1}}_{1} &=& \chi^1_1 & &
        \dot{\xi}^{k^{\xi}_{m}}_{m} &=& \chi^1_m \\[2mm]
      \end{array}\right.
    \end{equation}
    \begin{equation}\label{eq:TFCCS+IC_eta}
    \hspace{-2em}\Sigma_\chi:\left\{\hspace{1.5em}
      \begin{array}{rclcrclcrcl}
        \dot{\chi}^0 &=& \zeta^1_0 & &
        \dot{\chi}^1_1 &=& \chi^2_1 \zeta^1_0 + a^1_1(\xi,\chi^0,\bar\chi_2) & &
        \dot{\chi}^1_m &=& \chi^2_m \zeta^1_0 + a^1_m(\xi,\chi^0,\bar\chi_2) \\
        &&&& \dot{\chi}^2_1 &=& \chi^3_1 \zeta^1_0 + a^2_1(\xi,\chi^0,\bar\chi_3) & &
              \dot{\chi}^2_m &=& \chi^3_m \zeta^1_0 + a^2_m(\xi,\chi^0,\bar\chi_3) \\
        &&&& & \vdots & & \cdots & & \vdots & \\
        &&&& \dot{\chi}^{k^{\chi}-1}_1 &=& \chi^{k^{\chi}}_1 \zeta^1_0 + a^{k^{\chi}-1}_1(\xi,\chi) & &
              \dot{\chi}^{k^{\chi}-1}_m &=& \chi^{k^{\chi}}_m \zeta^1_0 + a^{k^{\chi}-1}_m(\xi,\chi) \\
        &&&& \dot{\chi}^{k^{\chi}}_1 &=& \zeta^1_1 + b^{k^\chi}_{0,1}(\xi, \chi)\zeta^1_0 & &
              \dot{\chi}^{k^{\chi}}_m &=& \zeta^1_m + b^{k^\chi}_{0,m}(\xi, \chi)\zeta^1_0 \\[2mm]
      \end{array}\right.
    \end{equation}
    \begin{equation}\label{eq:TFCCS+IC_zeta}
    \hspace{-7.7em}\Sigma_\zeta:\left\{ \hspace{0em}
      \begin{array}{rclcrclcrcl}
        \dot{\zeta}^1_0 &=& \zeta^2_0 & &
        \dot{\zeta}^1_1 &=& \zeta^2_1 & &
        \dot{\zeta}^1_m &=& \zeta^2_m \\
        & \vdots & & \hspace{0.2em} & & \vdots & \hspace{5.63em} & \cdots & \hspace{5.63em} & \vdots & \\
        \dot{\zeta}^{k^\zeta-s}_0 &=& w^0 & &
        \dot{\zeta}^{k^\zeta-s}_1 &=& \zeta^{k^\zeta-s+1}_1 & &
        \dot{\zeta}^{k^\zeta-s}_m &=& \zeta^{k^\zeta-s+1}_m \\
        & & & & & \vdots & & & & \vdots & \\
        & & & & \zetadot^{k^\zeta}_1 & = & w^1 & & \zetadot^{k^\zeta}_m & = & w^m \\[1mm]
      \end{array}\right.
    \end{equation}
\end{subequations}
\setcounter{equation}{\value{MYtempeqncnt}}
\end{figure*}

\section{Known Results}
\label{sec:known_results}
For the canonical contact form with compatible drift \eqref{eq:CCFD}, the drift and input vector fields can be written as\vspace{-0.5ex}
\begin{equation*}
\begin{aligned}
        a & = \sum_{i=1}^{k-1}\left(a^i_1\pad{z^i_1} + \ldots + a^i_m\pad{z^i_m} \right), \\ 
        b_0 & = \pad{z^0} + \sum_{i=1}^{k-1} \left( z^{i+1}_1\pad{z^{i}_1} + \ldots + z^{i+1}_m\pad{z^{i}_m} \right), b_j = \pad{z^k_j},\\
\end{aligned}\vspace{-0.5ex}
\end{equation*}
for $j=1,\ldots,m$. Necessary and sufficient conditions under which a control-affine system of the form \eqref{eq:ai_sys_m} is SFE to \eqref{eq:CCFD} are presented in \cite{LiMultiinputControlaffineSystems2016}. At generic points, these conditions can be stated as follows.
\begin{thm}
    Consider a control-affine system \eqref{eq:ai_sys_m} with $mk+1$ states, $m+1\ge3$ inputs, the drift $f$ and the input distribution $\Dist = \operatorname{span}\{g_0,\ldots,g_m\}$. The given system  is locally SFE to \eqref{eq:CCFD} if and only if:
    \begin{enumerate}[label=\arabic*.)]
        \item $\Dist^{(k-1)}=\TX$.
        \item $\Dist^{(k-2)}$ has rank $m(k-1)+1$ and contains an involutive subdistribution $\Ld$ of corank one.
        \item The drift $f$ satisfies the compatibility condition
        \begin{equation*}
            [f, \Cauchy{\Dist^{(i)}}] \subset \Dist^{(i)}, \quad \text{ for } \quad 1 \leq i \leq k-2 \; .
        \end{equation*}
    \end{enumerate}
\end{thm}

As shown in \cite{RespondekCanonicalContactSystems2002}, the existence of an involutive subdistribution $\Ld$ with corank one in $\Dist^{(k-2)}$ completely describes the geometry of systems SFE to \eqref{eq:CCFD} as summarized by the diagram
\begin{equation*}\label{diag:mcf}
    \arraycolsep=1pt
    \begin{array}{ccccccccccc}
        \Dist & \underset{m}{\subset} & \cdots & \underset{m}{\subset} & \Dist^{(k-3)} & \underset{m}{\subset} & \Dist^{(k-2)} & \underset{m}{\subset} & \Dist^{(k-1)} \hspace{-0.2ex} = \hspace{-0.2ex} \TX  \\[-0.2em] 
        \rotsubset &  & & & \rotsubset & & \rotsubset & & \\[0.9em]
        \Cauchy{\Dist^{(1)}} & \underset{m}{\subset} & \cdots & \underset{m}{\subset} & \Cauchy{\Dist^{(k-2)}} & \underset{m}{\subset} & \hspace{-1em}\Ld \; . &
    \end{array}
\end{equation*}

Since we aim to decompose a system \eqref{eq:ai_sys_m} into a contact part and integrator chains, we conclude this section with the following result from \cite{schlacher_jet_2015}, see also \cite[Lemma 2.2]{gstottner_analysis_2023}, to test whether a system allows a control-affine form.
\begin{lem}\label{lem:control_affine}
    Consider a system of the form \eqref{eq:gen_nonlin_sys} and define \mbox{$\Dist_0=\operatorname{span}\{\pad{u} \}$} and $\Dist_1=\Dist_0 + [f,\Dist_0]$. The system admits a control-affine representation if and only if \mbox{$\Dist_0 \subset \Cauchy{\Dist_1}$}.
\end{lem}
\vspace{-1em}
\section{Main Results}
\label{sec:main_results}
A natural augmentation of \eqref{eq:CCFD} is to prepend integrator chains to the inputs, i.e., to perform \emph{input prolongations}. A complete characterization of the fully general case, where arbitrary chain lengths at all inputs are allowed, appears computationally intractable. Therefore, we adopt a more tractable setting. We apply $k^\zeta$-fold input prolongations to each of the \(m\) chained inputs in \eqref{eq:CCFD}, except that the integrator chain associated to the distinguished input $w^0$ is shorter by $s$ (i.e., of length $k^\zeta-s$). In addition, we append integrator chains of arbitrary lengths to the top state variables. 
\addtocounter{equation}{1} The resulting system has the triangular structure \eqref{eq:TFCCS+IC_q}, with $\bar\chi_j=(\chi^1_1,\ldots,\chi^1_m,\ldots,\chi^j_1,\ldots,\chi^j_m)$ for $1\le j \le k^\chi$. We refer to \eqref{eq:TFCCS+IC_q} by \(\mathrm{TF}_s\). In general, this structurally flat triangular form consists of three subsystems:
\begin{itemize}
        \item \emph{The lower subsystem $\Sigma_\zeta$} is in Brunovský normal form. More precisely, it is given by $m$ integrator chains of length $k^\zeta$, and one integrator chain of length \mbox{$k^\zeta-s$}. Hence, $\Sigma_\zeta$ consists of $(m+1)k^\zeta-s$ states, with its top variables $\zeta^1_0, \ldots, \zeta^1_m$ acting as inputs for $\Sigma_\chi$.\footnote{In the case $s=k^\zeta$, the distinguished input $w^0$ enters $\Sigma_\chi$ directly, so there is no corresponding integrator chain in $\Sigma_\zeta$.}
        \item \emph{The middle system} $\Sigma_\chi$ is essentially of the form \eqref{eq:CCFD}, with the minor distinction that the functions $a^i_j$, for \mbox{$1\leq i \leq  k^\chi-1$} and $1 \leq j \leq m$, may also depend on the state $\xi$ of the upper subsystem $\Sigma_\xi$. We assume $k^\chi \geq 2$. This subsystem has $mk^\chi +1$ state variables and the top states $\chi^0, \chi^1_1 \ldots, \chi^1_m$ act as inputs for $\Sigma_\xi$. 
        \item  \emph{The upper subsystem} $\Sigma_\xi$ consists of $m+1$ integrator chains of arbitrary lengths $k^\xi_j\geq 0$, for $0 \leq j \leq m$, comprising $\sum_{j=0}^{m} k^\xi_j$ state variables.
    \end{itemize}
Hence, $\mathrm{TF}_s$ represents a system with $m+1$ inputs and \vspace{-1.5ex}
\begin{equation*}
    n = (m+1)k^\zeta-s + mk^\chi+1 + \sum_{j=0}^{m}k^\xi_j\vspace{-1.5ex}
\end{equation*}
states. 

In the following, we state necessary and sufficient for a control-affine system \eqref{eq:ai_sys_m} with $n \ge 2m+1$ states---we assume $k^\chi\ge2$---and $m+1 \ge 3$ inputs to be SFE to TF$_0$ and TF$_1$, respectively. The core idea is to characterize each subsystem of \eqref{eq:TFCCS+IC_q} individually and to capture their interconnection through supplementary conditions, thereby leading to a sequence of \emph{feedback invariant distributions}. We begin with the case where the lower subsystem \(\Sigma_{\zeta}\) consists of \(m+1\) integrator chains of equal length \(k^{\zeta}\), i.e. $s=0$. 
\begin{thm}\label{thm:TFCCS+IC_0}
    Consider a control-affine system \eqref{eq:ai_sys_m} with $m+1\ge3$ inputs. Let the drift $f$ and the distribution $\Dist_1=\Span{g_0,\ldots,g_m}$ define the sequence 
    \begin{equation}\label{eq:thm_TFCCS+IC_0_1_sequence}
        \quad \Dist_{i+1} = \Dist_i + [f,\Dist_i],\quad  i=1,\ldots, k^\zeta,
    \end{equation}
    with $k^\zeta$ being the smallest integer such that $\Dist_{k^\zeta+1}$ is non-involutive. The given system is SFE to TF$_0$ if and only if:
    \begin{enumerate}[label=\arabic*)]
        \item $\Rank{\Dist_i} = (m+1)i$, for all $1 \le i \le k^\zeta+1$.
        \item $\Cauchy{\Dist_{k^\zeta+1}}=\Dist_{k^\zeta}$.
    \end{enumerate}
    For the remaining conditions, let $\E = \Dist_{k^\zeta+1}$.
    \begin{enumerate}[resume]
        \item Let $k^\chi$ be the smallest integer such that $\E^{(k^\chi-1)} = \bar{\E}$. Then the following conditions must hold:
        \begin{enumerate}[label=\alph*)]
            \item $\operatorname{rank}(\bar{\E}) = (m+1)k^{\zeta} + mk^{\chi}+1.$
            \item $\operatorname{rank}(\E^{(k^\chi-2)}) =  \operatorname{rank}(\bar{\E})-m$ and there exists an involutive subdistribution $\Ld \crksub{1} \E^{(k^\chi - 2)}.$
            \item The drift $f$ satisfies the compatibility conditions:
            \begin{subequations}\label{eq:thm0_drift_comp}
                \begin{equation}\label{eq:thm0_drift_comp_1}
                    [f, \Cauchy{\E^{(i)}}] \subset \E^{(i)}, \text{ for } 1 \leq i \leq k^\chi-2,
                \end{equation}
                \begin{equation}\label{eq:thm0_drift_comp_2}
                    [f, \Ld] \subset \bar{\E}.
                \end{equation}
            \end{subequations}
        \end{enumerate}
        \item With $\F_0=\bar{\E}$, each distribution of the sequence
        \begin{equation*}
            \F_{i+1}=\F_{i} + [f, \F_{i}],\quad  i=0,\ldots, k^\xi-1,
        \end{equation*}
        where $k^\xi$ is the smallest integer such that $\F_{k^\xi} = \TX$, is involutive.
    \end{enumerate}
\end{thm}

In simple terms, 1) assures that $\Sigma_\zeta$ admits the Brunovský normal form and 2) that the top variables of $\Sigma_\zeta$ enter $\Sigma_\chi$ affinely. Conditions 3a), 3b) and \eqref{eq:thm0_drift_comp_1} ensure the contact form of $\Sigma_\chi$, while \eqref{eq:thm0_drift_comp_2} assures that the top variables of $\Sigma_\chi$ act as inputs for $\Sigma_\xi$. 
Finally, item~4) guarantees that $\Sigma_\xi$ can be transformed to integrator chains.
A sketch of the proof of Theorem~\ref{thm:TFCCS+IC_0} is given in Appendix \ref{app:proof_thm_TFCCS_0}. To verify the existence of an involutive corank one subdistribution $\Ld \crksub{1} \E^{(k^\chi-2)}$ required by condition~3b) of Theorem~\ref{thm:TFCCS+IC_0}, we follow the constructive procedure as described, e.g., in \cite{respondek_transforming_2001, RespondekCanonicalContactSystems2002, LiMultiinputControlaffineSystems2016} (for more details about involutive corank one subdistributions see \cite{pasillas-lepine_contact_2001}). The following Lemma summarizes this procedure.

\begin{lem}\label{lem:constructing_L}
Let $\Dist$ be a rank-$d$ distribution on the $n$-dimensional manifold $\X$ with its annihilator $\mathcal{P}$. For any $\omega^i \in \mathcal{P}$, define
\begin{equation*}
    \mathcal{W}_i=\{ v \in \Dist: v \hook \D \omega^i \in \Dist^\perp  \}=\Cauchy{\operatorname{span}\{\D \omega^i\}^\perp} \cap \Dist .
\end{equation*}
Assume $\operatorname{rank}(\Dist^{(1)}) = d + r$ (i.e., $\operatorname{rank}(\mathcal{P}^{(1)}) = n-d-r$) with $r\ge2$. Choose one-forms \mbox{$\omega^1, \ldots, \omega^r, \omega^{r+1}, \ldots, \omega^{n-d}$} such that \mbox{$\mathcal{P}=\operatorname{span}\{\omega^1, \ldots,  \omega^{n-d}\}$} and \mbox{$\mathcal{P}^{(1)} = \Span{\omega^{r+1},\ldots,\omega^{n-d}}$}. Next, build the sum of two distributions $\Ld = \mathcal{W}_i + \mathcal{W}_j \text{ for any } 1 \leq i < j \leq r \; $. If $\Ld$ is involutive and has corank one in $\Dist$, then it is the unique involutive corank one subdistribution of $\Dist$. If the resulting $\Ld$ does not have corank one in $\Dist$ or is not involutive, there does not exist any involutive distribution which has corank one in $\Dist$.
\end{lem}
\vspace{1ex}

Starting from $\mathrm{TF}_0$, the next more complicated case is given by \eqref{eq:TFCCS+IC_q} with the first integrator chain being shorter by $s=1$.
\begin{thm}\label{thm:TFCCS+IC_1}
    Consider a system \eqref{eq:ai_sys_m} with the drift $f$, the distribution $\Dist_1=\operatorname{span}\{g_0,\ldots,g_m\}$ and the sequence \eqref{eq:thm_TFCCS+IC_0_1_sequence}. The given system is SFE to $\mathrm{TF}_1$ if and only if:
    \begin{enumerate}[label=\arabic*)]
        \item $\operatorname{rank}(\Dist_{i}) = (m+1)i,\text{ for all } 1 \leq i \leq k^\zeta+1$.
        \item $\Cauchy{\Dist_{k^\zeta+1}}\neq\Dist_{k^\zeta}$.
        \item There exist $m$ vector fields \mbox{$c_i\in \Dist_{k^\zeta}$}, \mbox{$c_i\notin \Dist_{k^\zeta-1}$}, \mbox{$i=1,\ldots,m$}, which are linearly independent modulo $\Dist_{k^\zeta-1}$, such that
        \begin{equation}\label{eq:thm_TFCCS+IC_1_item_1_E}
            \E = D_{k^{\zeta}} + \operatorname{span}\{[f, c_1], \ldots, [f, c_{m}]\}    
        \end{equation}
        satisfies
        \begin{equation}\label{eq:thm_EMCF+IC_1_item_1_CE}
            \Cauchy{\E} = D_{k^\zeta - 1} + \operatorname{span}\{c_1,\ldots,c_{m}\}\; .
        \end{equation}
        \item Condition~3) of Theorem~\ref{thm:TFCCS+IC_0} holds for \eqref{eq:thm_TFCCS+IC_1_item_1_E} with the difference that $\operatorname{rank}(\bar{\E}) = (m+1)k^{\zeta} + mk^{\chi}$.
        \item With $\F_0=\bar \E$, condition~4 as defined in Theorem~\ref{thm:TFCCS+IC_0} holds.
    \end{enumerate} 
\end{thm}
\vspace{1ex}

Within Theorem~\ref{thm:TFCCS+IC_1}, $\Sigma_\zeta$ is characterized via the involutive distributions \mbox{$\Dist_1 \subset \cdots \Dist_{k^\zeta-1} \subset \Cauchy{\E}$}. By constructing the corresponding vector fields $c_1,\ldots,c_m$, we identify the $m$ integrator chains whose top variables act as chained inputs for $\Sigma_\chi$. The subsystems $\Sigma_\chi$ and $\Sigma_\xi$ are characterized the same way as in Theorem~\ref{thm:TFCCS+IC_0}. For the outline of the proof, see Appendix~\ref{app:proof_thm_TFCCS_1}. Condition~3 of Theorem~\ref{thm:TFCCS+IC_1} requires the construction of $m$ specific vector fields $c_1,\ldots,c_{m}$. For $m=1$, \cite[Sec. 4.1.1]{gstottner_structurally_2022} provides necessary conditions from which candidates for $c_1$ are derived purely by algebraic operations and differentiation. The next lemma extends these conditions to arbitrary $m \ge 1$.
\begin{lem}\label{lem:TFCCS_1}
    Consider a sequence of distributions \mbox{$\Dist_0 \underset{m+1}{\subset} \Dist_1 \underset{m+1}{\subset} \Dist_2$} where $\Dist_1$ is involutive and $\Dist_2$ is non-involutive. Assume further that there exists a vector field $f$ such that $[f, \Dist_0] \subset \Dist_1$ and $\Dist_2 = \Dist_1 + [f, \Dist_1]$ and suppose that $\Dist_1 = \Dist_0 + \Span{v_0, \ldots, v_m}$. Consider the vector fields
    \begin{equation}\label{eq:tech_result_vec_fields_c}
        c_i = \sum^{m}_{k=0} \alpha^k_iv_k , \quad i=1,\ldots,m,
    \end{equation}
    with coefficients $\alpha^k_i$, such that \eqref{eq:tech_result_vec_fields_c} are linearly independent modulo $\Dist_0$. Let \mbox{$\E_1 = \Dist_0 + \Span{c_1, \ldots, c_{m}}$} and \mbox{$\E_2 = \Dist_1 + \Span{[c_1,f], \ldots, [c_{m}, f]}$}. If there exist $m$ vector fields of the form \eqref{eq:tech_result_vec_fields_c} such that $\E_1 \subset \Cauchy{\E_2}$, then the associated $(m+1)m$ coefficient functions $\alpha^k_i$ must satisfy
    \begin{equation}\label{eq:tech_result_conditions_set}
        \sum_{k=0}^{m} \alpha^k_i \alpha^k_j [v_k,[v_k,f]] + \hspace{-1.1em} \sum_{0\le l < p \le m} \hspace{-0.9em} (\alpha ^l_i \alpha^p_j + \alpha^p_i \alpha ^l_j)\,[v_l,[v_p,f]]
        \hspace{-0.2em} \in \hspace{-0.2em} \Dist_2
    \end{equation}
    for $1 \le i \le j \le m$, yielding $\binom{m+1}{2}$ conditions in total.
\end{lem}
\vspace{1ex}

The proof of Lemma \ref{lem:TFCCS_1} is omitted for brevity but can be conducted analogously to the argumentation of \cite[Sec. 4.1.1]{gstottner_structurally_2022}. By applying Lemma \ref{lem:TFCCS_1} to the sequence \mbox{$\Dist_{k^\zeta-1} \crksub{m+1} \Dist_{k^\zeta} \crksub{m+1} \Dist_{k^\zeta+1}$} and solving \eqref{eq:tech_result_conditions_set} for the \mbox{$(m+1)m$} unknown functions $\alpha^k_i$, one obtains candidates for the vector fields $c_i$, $i=1,\ldots,m$, in item~3) of Theorem~\ref{thm:TFCCS+IC_0}.

\subsection{Computation of Flat Outputs}

Systems in TF\(_0\) or TF\(_1\) are flat, with the top variables of the triangular form providing a compatible flat output.  If a system \eqref{eq:ai_sys_m} is SFE to TF\(_0\) or TF\(_1\), explicitly transforming the system into TF$_0$ or TF$_1$ is not required to compute such flat outputs. All compatible flat outputs, i.e., all sets of functions which can serve as top variables of the corresponding triangular form, can be obtained directly from the sequence of involutive distributions \mbox{$\Ld \subset \F_0 \subset \cdots \subset \F_{k^\xi-1}$} appearing in the Theorems~\ref{thm:TFCCS+IC_0} and~\ref{thm:TFCCS+IC_1}. 

\section{A Practical Example}
\label{sec:example}
\begin{figure}[!t]
    \centering
    \def\svgwidth{0.45\textwidth}
    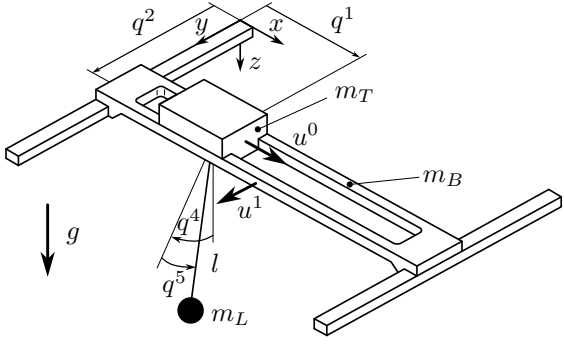
    \caption{ Schematic diagram of a gantry crane. }
    \label{fig:crane}
\end{figure}

This section demonstrates static feedback equivalence of  the three-dimensional gantry crane to $\mathrm{TF}_1$ by verifying the conditions of Theorem~\ref{thm:TFCCS+IC_1}, thereby, applying Lemmas \ref{lem:constructing_L} and \ref{lem:TFCCS_1}, and transforming the crane into  $\mathrm{TF}_1$.
The gantry crane as shown in Figure \ref{fig:crane} is a holonomic mechanical system with five degrees of freedom and three control inputs. The coordinates $q^1$ and $q^2$ denote the position of the trolley in the $xy$-plane. The length of the rope is given by $l=R_dq^3$ with the radius of the rope drum $R_d$ and its rotation angle $q^3$. Furthermore, let $q^4$ be the rope’s swing angle in the \(yz\)-plane, and $q^5$ the angle between the rope and its \(yz\)-projection. Accordingly, the position of the load is given by
\begin{equation}\label{eq:load_pos}
    \begin{aligned}
        x_L = & q^1 + R_dq^3\sin(q^5),\\
        y_L = & q^2 + R_dq^3\sin(q^4)\cos(q^5), \\
        z_L = & R_dq^3\cos(q^4)\cos(q^5).
    \end{aligned}
\end{equation}

To derive the state-space model of the crane, we neglect friction and apply the Lagrangian formalism using the Lagrangian \mbox{$L=T-V$}. With $m_T$ and $m_B$ indicating the mass of the trolley and the bridge and $J$ denoting the inertia of the drum, $T$ refers to the kinetic energy given by
\begin{equation*}
\resizebox{\linewidth}{!}{$
    T = \frac{m_L}{2}(\dot{x}_L^2+\dot{y}_L^2+\dot{z}_L^2)+\frac{J}{2}(v^3)^2+\frac{m_B+m_T}{2}(v^2)^2+\frac{m_T}{2}(v^1)^2. 
    $}
\end{equation*}
Further,  $V = -m_Lgz_L$ stands for the potential energy. with the gravitational acceleration $g$. The velocities corresponding to the configuration coordinates are denoted by  $v=(v^1,\ldots,v^5)$ and complement $q$ to the full state $x=(q,v)$. 

Consequently, the crane's control-affine state-space representation is of the form
\begin{equation}\label{eq:gantry_crane_sys}
        \dot{q} = v, \hspace{1em} \dot{v} = f_v + g_0u^0 + g_1u^1 + g_2u^2 ,
\end{equation}
with the forces $u^0$ and $u^1$ and the torque $u^2$ acting on the rope drum. The drift and the input vector fields of \eqref{eq:gantry_crane_sys} are of the form
\begin{subequations}
    \begin{alignat}{1}
        f & = v^i\pad{q^i} + f_v^i(q^3, q^4, q^5, v^4, v^5)\pad{v^i}, \label{eq:drift_crane} \\
        g_j & = g^i_j(q^3,q^4,q^5)\pad{v^i}, \quad  \label{eq:input_crane}
    \end{alignat}
\end{subequations}
with $i=1,\ldots,5$, $j=0,\ldots,2$. Next, we show SFE of \eqref{eq:gantry_crane_sys} to TF$_1$ by verifying the conditions of Theorem~\ref{thm:TFCCS+IC_1}.

\emph{Computing the sequence \eqref{eq:thm_TFCCS+IC_0_1_sequence}.} Given the involutive input distribution $\Dist_1=\Span{g_0,g_1,g_2}$, we compute \mbox{$\Dist_2 = \Dist_1 + [f,\Dist_1]$}. Since $\Dist_2$ is already non-involutive, we obtain $k^\zeta=1$. 

\emph{Condition~1) and 2).} Given $\Rank{\Dist_2}=6$ and \mbox{$\Cauchy{\Dist_2} \neq \Dist_1$}, item~1) and 2) of Theorem~\ref{thm:TFCCS+IC_1} are satisfied. 

\emph{Condition~3).} To construct $m=2$ vector fields $c_1$ and $c_2$ such that \eqref{eq:thm_TFCCS+IC_1_item_1_E} satisfies \eqref{eq:thm_EMCF+IC_1_item_1_CE}, we apply Lemma \ref{lem:TFCCS_1}. For the gantry crane, the sequence of distributions  $\{ \emptyset \} \underset{3}{\subset} \Dist_1 \underset{3}{\subset} \Dist_2$ meets the assumptions of Lemma \ref{lem:TFCCS_1} with the drift vector field \eqref{eq:drift_crane}. Accordingly, we define the vector fields \eqref{eq:tech_result_vec_fields_c} as
\begin{equation}\label{eq:crane_c_i}
    c_i = \alpha^0_ig_0 + \alpha^1_ig_1 + \alpha^2_ig_2, \quad i=1,2.
\end{equation}
With the solution
\begin{equation}\label{eq:crane_alpha_1}
    \begin{aligned}
        \alpha^0_1 &= -R_dq^3\cos(q^4)\cos(q^5)m_T, \; \alpha^1_1 = 0, \\
        \; \alpha^2_1 &= \frac12 Jq^3\sin(2q^5)\cos(q^4) 
    \end{aligned}
\end{equation}
and
\begin{equation}\label{eq:crane_alpha_2}
    \begin{aligned}
        \alpha^0_2 &= 0,\; \alpha^1_2 = -(m_B+m_T)Rq^3\cos(q^5)\cos(q^4),\\
        \alpha^2_2 &= \frac12Jq^3\sin(2q^4)\cos^2(q^5)
    \end{aligned}
\end{equation}
of the three corresponding algebraic conditions \eqref{eq:tech_result_conditions_set}, we obtain vector fields $c_1, c_2$ such that the distribution
\begin{equation*}
    \E = \Span{g_0, g_1, g_2, [f,c_1], [f,c_2]}
\end{equation*}
satisfies \mbox{$\Cauchy{\E} = \Span{c_1, c_2}$}. Any other candidates for the vector fields $c_1,c_2$ obtained via solving \eqref{eq:tech_result_conditions_set} yield the same distribution $\E$.

\emph{Condition~4).} With $\E^{(1)}=\bar{\E}$, it follows that $k^\chi=2$. Given that $\Rank{\bar\E} = 7=(m+1)k^\zeta + mk^\chi$, it remains to verify condition~3b) and 3c) of Theorem~\ref{thm:TFCCS+IC_0}. The rank condition of item~3b) is met with $\Rank{\E} = 5$. To proceed, we aim to construct an involutive subdistribution $\Ld \underset{1}{\subset} \E$, such that item~3b) is satisfied entirely.
By choosing two linearly independent one-forms $\omega^1,\omega^2$ such that $\omega^i \in \E^\perp, \omega^i\notin (\E^{(1)})^\perp, i=1,2,$ we construct the distributions $\mathcal{W}_1$ and $\mathcal{W}_2$. By building the union of these two distributions, we obtain the subdistribution $\Ld = \mathcal{W}_1 + \mathcal{W}_2$ that is involutive and has corank one in $\E$, i.e., $\Rank{\Ld} = 4$.
For the remaining item~3c), we have to check the compatibility conditions of the drift vector field. This is accomplished by showing that $[f,\Ld] \subset \bar{\E}$.

\emph{Condition~5).} With $\F_0 = \bar{\E}$, we obtain $\F_1= \TX$. Consequently, the gantry crane meets also item~5) of Theorem~\ref{thm:TFCCS+IC_1}. Note that $\Rank{\F_1} - \Rank{\F_0} = 3$ implies that $\Sigma_\xi$ consists of 3 integrator chains of length $k^\xi=1$.

In fact, flat outputs of the gantry crane \eqref{eq:gantry_crane_sys} are all triples \(\varphi=(\varphi^0,\varphi^1,\varphi^2)\) whose differentials span~\(\F_0^\perp\), i.e., $\operatorname{span}\{\D\varphi^0,\D\varphi^1,\D\varphi^2\}=\F_0^\perp$.
In the coordinates of the state-space representation \eqref{eq:gantry_crane_sys}, the annihilator of the involutive distribution $\F_0$ is given by
\begin{equation*}
\begin{aligned}
    \F_0^\perp = & \operatorname{span} \left\{ \frac{\cos{q^5}}{R_dq^3}\D q^1 - \frac{\sin(q^4)\sin(q^5)}{R_dq^3}\D q^2 + \D q^5, \right. \\
    & \frac{\cos(q^4)}{R_dq^3\cos(q^5)}\D q^2 + \D q^4,\\
    & \left. \frac{\sin(q^5)}{R_d}\D q^1 + \frac{\sin(q^4)\cos(q^5)}{R_d}\D q^2 + \D q^3 \right\}.
\end{aligned}
\end{equation*}

In particular, $\F_0^\perp$ is spanned by the differentials of \eqref{eq:load_pos}. Hence, the position of the load is a TF$_1$-compatible flat output.
To conclude, we have systematically derived a flat output that is compatible with $\mathrm{TF}_1$ by means of the geometric characterization provided in Theorem~\ref{thm:TFCCS+IC_1}. This complements the previous equation-driven derivation of the same flat output in \cite{rudolph_examples_1998}, which analyzes a reduced crane model and offers an insightful, example-oriented perspective.

Having shown that the gantry crane is SFE to TF\(_1\), we now carry out the state and input transformations that bring the system into TF\(_1\), following the sufficiency steps of the sketch of the proof for Theorem~\ref{thm:TFCCS+IC_1}.

\emph{Step~1 \& 2: Decomposing \eqref{eq:gantry_crane_sys} into $\Sigma_\zeta, \Sigma_\chi, \Sigma_\xi$}. By means of the vector fields \eqref{eq:crane_c_i} with its corresponding coefficient functions \eqref{eq:crane_alpha_1} and \eqref{eq:crane_alpha_2}, we introduce the input transformation \mbox{$(u^0, u^1, u^2)=$} \mbox{$ (\alpha^0_1\bar{u}^1 + \alpha^0_2\bar{u}^2, \alpha^1_1\bar{u}^1 + \alpha^1_2\bar{u}^2, \alpha^2_1\bar{u}^1 + \alpha^2_2\bar{u}^2 + \bar{u}^0 )$}, such that \mbox{$\xdot = f(x) + g_2(x)\bar{u}^0 + c_1(x)\bar{u}^1 + c_2(x)\bar{u}^2$}. We chose this approach, to split off the inputs $(\bar{u}^1, \bar{u}^2)$ when straightening out the vector fields $c_1,c_2$. 
From the geometric characterization of Theorem~\ref{thm:TFCCS+IC_1}, we obtain the sequence of involutive distributions
\begin{equation*}\label{eq:seq_dist_gantry}
    \Cauchy{\E} \underset{2}{\subset} \Ld \underset{3}{\subset} \bar{\E} = \F_0 \underset{3}{\subset} \F_1 = \TX \; ,
\end{equation*}
where $\Ld^\perp = \operatorname{span}\{ \D \varphi^0, \D \varphi^1, \D \varphi^2, \D \Lie_{f} \varphi^0, \D \Lie_{f} \varphi^1, \D \Lie_{f} \varphi^2 \}$. 
To straighten out these distributions simultaneously, we apply the state transformation
\begin{equation*}
    \arraycolsep = 1.4pt
    \begin{array}{rclcrclcrcl}
         \xi^1_0 & = & \varphi^2, & & \xi^1_1 & = & \varphi^0, & & \xi^1_2 & = & \varphi^1, \\[0.5ex]
         \bar{\chi}_0^1 & = & \Lie_f \varphi^2, & & \bar{\chi}^1_1 & = & \Lie_f \varphi^0, & & \bar{\chi}^1_2 & = & \Lie_f \varphi^1, \\[0.5ex]
         & & & & \bar{\chi}^2_1 & = & q^4, & & \bar{\chi}^2_2 & = & q^5, \\[0.5ex]
         & & & & \bar{\zeta}^1_1 & = & v^4, & & \bar{\zeta}^1_2 & = & v^5 . \\[0.5ex]
    \end{array}
\end{equation*}
Consequently, we have decomposed the system into three subsystems $\Sigma_\xi, \Sigma_\chi,\Sigma_\zeta$ of the form
\begin{equation*}
    \arraycolsep = 1.4pt
    \begin{array}{rclcrclcrcl}
        \xidot^1_0 & = & \bar{\chi}^1_0, & & \xidot^1_1 & = & \bar{\chi}^1_1, & & \xidot^1_2 & = & \bar{\chi}^1_2, \\[0.5ex]
        \dot{\bar{\chi}}^1_0 & = & f^1_{\bar \chi,0}(\tilde{z}_0, \bar{u}^0), & & \dot{\bar{\chi}}^1_1 & = & f^1_{\bar \chi,1}(\tilde{z}_0, \bar{u}^0), & & \dot{\bar{\chi}}^1_2 & = & f^1_{\bar \chi,2}(\tilde{z}_0, \bar{u}^0), \\[0.5ex]
        & & & & \dot{\bar{\chi}}^2_1 & = & \bar{\zeta}^1_1, & & \dot{\bar{\chi}}^2_2 & = & \bar{\zeta}^1_2, \\[0.5ex]
        & & & & \dot{\bar{\zeta}}^1_1 & = & f^1_{\bar \zeta,1}(\tilde{z}_1, \bar{u}^0, \bar{u}^2), & & \dot{\bar{\zeta}}^1_2 & = & f^1_{\bar \zeta,2}(\tilde{z}_1, \bar{u}),  \\
    \end{array}
\end{equation*}
where $\tilde{z}_0 \hspace{-0.5em}=\hspace{-0.5em} (\xi^1_0, \bar{\chi}_1^2, \bar{\chi}_2^2, \bar{\zeta}^1_1, \bar{\zeta}^1_2)$ and $\tilde{z}_1 \hspace{-0.5em}= \hspace{-0.5em}(\xi^1_0, \bar{\chi}_0^1, \bar{\chi}_1^2, \bar{\chi}_2^2, \bar{\zeta}^1_1, \bar{\zeta}^1_2)$. Given that $\Sigma_\xi$ is already in Brunovský normal form, we proceed with step~4 of the sufficiency part of Appendix \ref{app:proof_thm_TFCCS_1}. 

\emph{Step~4: Transforming $\Sigma_\chi$ to the contact form.}  
By applying the input transformation $w^0 = f^1_{\bar \chi,0}(\tilde{z}_0, \bar{u}^0)$ with the inverse $\bar{u}^0 = \Phi_{\bar u^0}(\tilde{z}_0,w^0)$, we render $\Sigma_{\chi}$ into a control affine system with the inputs $(w^0, \bar{\zeta}^1_1, \bar{\zeta}^1_2)$:

\begin{equation*}
    \arraycolsep = 1.4pt
    \begin{array}{rclcrclcrcl}
         \dot{\bar{\chi}}^1_0 & = & w^0, & & \dot{\bar{\chi}}^1_1 & = & (w^0-g)\frac{\tan(\bar{\chi}^2_2)}{\cos(\bar{\chi}^2_1)}, & & \dot{\bar{\chi}}^1_2 & = & (w^0-g)\tan(\bar{\chi}^2_1), \\[0.5ex]
         &&&& \dot{\bar{\chi}}^2_1 & = & \bar{\zeta}^1_1, & & \dot{\bar{\chi}}^2_2 & = & \bar{\zeta}^1_2  .\\
    \end{array}
\end{equation*}
By applying the state transformation $(\chi^0, \chi^1_1, \chi^1_2, \chi^2_1, \chi^2_2) = (\bar \chi_0^1, \bar \chi^1_1, \bar \chi^1_2,\frac{\tan(\bar{\chi}^2_2)}{\cos(\bar{\chi}^2_1)}, \tan(\bar{\chi}^2_1))$ we obtain
\begin{equation*}
    \arraycolsep = 1.4pt
    \begin{array}{rclcrclcrcl}
         \dot{\chi}^0 & = & w^0, & & \dot{\chi}^1_1 & = & (w^0-g)\chi^2_1, & & \dot{\chi}^1_2 & = & (w^0-g)\chi^2_2, \\[0.5ex]
         &&&& \dot{\chi}^2_1 & = & f^2_{\chi,1}(\chi_1^2, \chi_2^2, \bar \zeta^1_1, \bar \zeta^1_2 ), & & \dot{\chi}^2_2 & = & f^2_{\chi,2}(\chi^2_2, \bar \zeta^1_1)  .\\
    \end{array}
\end{equation*}
Finally, by applying the transformation \mbox{$(\zeta^1_1, \zeta^1_2) = (f^2_{\chi,1},f^2_{\chi,2})$}, subsystem $\Sigma_\chi$ is of the form \eqref{eq:TFCCS+IC_eta} with $k^\chi=2$.

\emph{Step~5: Transforming $\Sigma_\zeta$ into Brunovský normal form.} In the last step, $\Sigma_\zeta$, given by the state equations $\dot \zeta^1_1 = f^1_{\zeta,1}( \xi^1_0, \chi^1_0, \chi^2_1, \zeta^1_1, w^0, \bar{u}^1 )$ and $\dot \zeta^1_2 = f^1_{\zeta,2}( \xi^1_0, \chi^1_0, \chi^2_2, \zeta^1_2, w^0, \bar{u}^2 )$, is transformed into two integrator chains of length $k^\zeta=1$ by simply applying the input transformation $(w^1, w^2) = (f^1_{\zeta,1}, f^1_{\zeta,2})$. 

Finally, the gantry crane in $\mathrm{TF}_1$ is given by
\begin{alignat*}{3}
    \xidot^1_0   &= \chi^0,    &\hspace{3em} \xidot^1_1   &= \chi^1_1,    &\hspace{3em} \xidot^1_2   &= \chi^1_2,   \\[1ex]
    \chidot^0  &= w^0,         &\quad \chidot^1_1  &= \chi^2_1 w^0 - \chi^2_1 g, 
                                &\quad \chidot^1_2  &= \chi^2_2 w^0 - g \chi^2_2,   \\
                 &             &\quad \chidot^2_1  &= \zeta^1_1, 
                                &\quad \chidot^2_2  &= \zeta^1_2,  \notag \\[1ex]
    & &\quad \zetadot^1_1 &= w^1,         &\quad \zetadot^1_2 &= w^2  .
\end{alignat*}

\section{Conclusion}
\label{sec:conclusion}

We introduced the triangular form TF$_s$ for multi-input control-affine systems and gave necessary and sufficient conditions for static feedback equivalence to TF$_0$ and TF$_1$. The conditions are constructive in the sense that compatible flat outputs can be computed directly from the involved distributions, without explicitly transforming the system. We demonstrated that the gantry crane in three dimensions is SFE to TF\(_1\), with the load position as flat output recovered systematically. Our results extend structurally flat triangular forms, and provide sufficient conditions for flatness as well as procedures to derive flat outputs compatible with TF\(_0\) and TF\(_1\).

\bibliography{IEEEabrv,mybibfile}

\begin{thebibliography}{10}
\providecommand{\url}[1]{#1}
\csname url@samestyle\endcsname
\providecommand{\newblock}{\relax}
\providecommand{\bibinfo}[2]{#2}
\providecommand{\BIBentrySTDinterwordspacing}{\spaceskip=0pt\relax}
\providecommand{\BIBentryALTinterwordstretchfactor}{4}
\providecommand{\BIBentryALTinterwordspacing}{\spaceskip=\fontdimen2\font plus
\BIBentryALTinterwordstretchfactor\fontdimen3\font minus \fontdimen4\font\relax}
\providecommand{\BIBforeignlanguage}[2]{{%
\expandafter\ifx\csname l@#1\endcsname\relax
\typeout{** WARNING: IEEEtran.bst: No hyphenation pattern has been}%
\typeout{** loaded for the language `#1'. Using the pattern for}%
\typeout{** the default language instead.}%
\else
\language=\csname l@#1\endcsname
\fi
#2}}
\providecommand{\BIBdecl}{\relax}
\BIBdecl

\bibitem{fliess_flatness_1995}
M.~Fliess, J.~Lévine, P.~Martin, and P.~Rouchon, ``\BIBforeignlanguage{en}{Flatness and defect of non-linear systems: introductory theory and examples},'' \emph{\BIBforeignlanguage{en}{International Journal of Control}}, vol.~61, no.~6, pp. 1327--1361, 1995.

\bibitem{fliess_lie-backlund_1999}
------, ``\BIBforeignlanguage{en}{A {Lie}-{Bäcklund} approach to equivalence and flatness of nonlinear systems},'' \emph{\BIBforeignlanguage{en}{IEEE Trans. Autom. Control}}, vol.~44, no.~5, pp. 922--937, May 1999.

\bibitem{gstottner_tracking_2024}
C.~Gstöttner, B.~Kolar, and M.~Schöberl, ``\BIBforeignlanguage{en}{Tracking control for (x, u)-flat systems by quasi-static feedback of classical states},'' \emph{\BIBforeignlanguage{en}{Symmetry, Integrability and Geometry: Methods and Applications (SIGMA)}}, vol.~20, no.~71, Jul. 2024.

\bibitem{schoberl_implicit_2014}
M.~Schöberl and K.~Schlacher, ``\BIBforeignlanguage{en}{On an implicit triangular decomposition of nonlinear control systems that are 1-flat—{A} constructive approach},'' \emph{\BIBforeignlanguage{en}{Automatica}}, vol.~50, no.~6, pp. 1649--1655, Jun. 2014.

\bibitem{nicolau_flatness_2017}
F.~Nicolau and W.~Respondek, ``Flatness of multi-input control-affine systems linearizable via one-fold prolongation,'' \emph{SIAM Journal on Control and Optimization}, vol.~55, no.~5, pp. 3171--3203, Jan. 2017.

\bibitem{gstottner_finite_2021}
C.~Gstöttner, B.~Kolar, and M.~Schöberl, ``\BIBforeignlanguage{en}{A finite test for the linearizability of two-input systems by a two-dimensional endogenous dynamic feedback},'' in \emph{\BIBforeignlanguage{en}{2021 {European} {Control} {Conference} ({ECC})}}, Delft, Netherlands, 2021, pp. 970--977.

\bibitem{gstottner_necessary_2023}
------, ``\BIBforeignlanguage{en}{Necessary and sufficient conditions for the linearisability of two-input systems by a two-dimensional endogenous dynamic feedback},'' \emph{\BIBforeignlanguage{en}{International Journal of Control}}, vol.~96, no.~3, pp. 800--821, 2023.

\bibitem{levine_differential_2025}
J.~Lévine, ``Differential flatness by pure prolongation: necessary and sufficient conditions,'' \emph{Communications in Optimization Theory}, vol. 2025, no.~11, pp. 1--40, 2025.

\bibitem{nicolau_dynamic_2025}
F.~Nicolau, W.~Respondek, and S.~Li, ``\BIBforeignlanguage{en}{Dynamic feedback linearization of two-input control systems via successive one-fold prolongations},'' \emph{\BIBforeignlanguage{en}{Journal of the Franklin Institute}}, vol. 362, no.~13, p. 107780, Aug. 2025.

\bibitem{martin_feedback_1994}
P.~Martin and P.~Rouchon, ``Feedback linearization and driftless systems,'' \emph{Mathematics of Control, Signals, and Systems}, vol.~7, no.~3, pp. 235--254, Sep. 1994.

\bibitem{murray_nilpotent_1994}
R.~M. Murray, ``\BIBforeignlanguage{en}{Nilpotent bases for a class of nonintegrable distributions with applications to trajectory generation for nonholonomic systems},'' \emph{\BIBforeignlanguage{en}{Mathematics of Control, Signals, and Systems}}, vol.~7, no.~1, pp. 58--75, Mar. 1994.

\bibitem{li_characterization_2013}
S.~Li, C.~Xu, H.~Su, and J.~Chu, ``Characterization and flatness of the extended chained system,'' in \emph{Proceedings of the 32nd {Chinese} {Control} {Conference}}, Xi'an, China, 2013.

\bibitem{bououden_triangular_2011}
S.~Bououden, D.~Boutat, G.~Zheng, J.-P. Barbot, and F.~Kratz, ``A triangular canonical form for a class of 0-flat nonlinear systems,'' \emph{International Journal of Control}, vol.~84, no.~2, pp. 261--269, 2011.

\bibitem{gstottner_flat_2021}
C.~Gstöttner, B.~Kolar, and M.~Schöberl, ``On a flat triangular form based on the extended chained form,'' \emph{IFAC-PapersOnLine}, vol.~54, no.~9, pp. 245--252, 2021.

\bibitem{gstottner_structurally_2022}
------, ``A structurally flat triangular form based on the extended chained form,'' \emph{International Journal of Control}, vol.~95, no.~5, pp. 1144--1163, 2022.

\bibitem{gstottner_triangular_2024}
------, ``\BIBforeignlanguage{en}{A triangular normal form for x-flat control-affine two-input systems},'' in \emph{\BIBforeignlanguage{en}{2024 28th {International} {Conference} on {Methods} and {Models} in {Automation} and {Robotics} ({MMAR})}}, Miedzyzdroje, Poland, Aug. 2024, pp. 298--303.

\bibitem{hartl_triangular_2025}
\BIBentryALTinterwordspacing
G.~Hartl, C.~Gstöttner, and M.~Schöberl, ``\BIBforeignlanguage{en}{On triangular forms for x-flat control-affine systems with two inputs},'' in \emph{\BIBforeignlanguage{en}{29th {International} {Conference} on {System} {Theory}, {Control} and {Computing} ({ICSTCC})}}, Cluj-Napoca, Romania, 2025, to appear. [Online]. Available: \url{https://arxiv.org/abs/2505.15562}
\BIBentrySTDinterwordspacing

\bibitem{murray_nonholonomic_1993}
R.~Murray and S.~Sastry, ``\BIBforeignlanguage{en}{Nonholonomic motion planning: steering using sinusoids},'' \emph{\BIBforeignlanguage{en}{IEEE Transactions on Automatic Control}}, vol.~38, no.~5, pp. 700--716, May 1993.

\bibitem{respondek_transforming_2001}
W.~Respondek, ``Transforming nonholonomic control systems into the canonical contact form,'' in \emph{Proceedings of the 40th {IEEE} {Conference} on {Decision} and {Control}}, vol.~2, Orlando, FL, USA, 2001, pp. 1781--1786.

\bibitem{RespondekCanonicalContactSystems2002}
W.~Respondek and W.~Pasillas-Lépine, ``\BIBforeignlanguage{en}{Canonical contact systems for curves: a survey},'' \emph{\BIBforeignlanguage{en}{Contemporary Trends in Nonlinear Geometric Control Theory and Its Applications}}, pp. 77--112, Jan. 2002.

\bibitem{LiMultiinputControlaffineSystems2016}
S.~Li, F.~Nicolau, and W.~Respondek, ``Multi-input control-affine systems static feedback equivalent to a triangular form and their flatness,'' \emph{International Journal of Control}, vol.~89, no.~1, pp. 1--24, 2016.

\bibitem{schlacher_jet_2015}
K.~Schlacher, M.~Schöberl, and B.~Kolar, ``\BIBforeignlanguage{en}{A {Jet} {Space} {Approach} to {Derive} {Flat} {Outputs}},'' \emph{\BIBforeignlanguage{en}{IFAC-PapersOnLine}}, vol.~48, no.~11, pp. 131--136, 2015.

\bibitem{gstottner_analysis_2023}
C.~Gstöttner, \emph{\BIBforeignlanguage{en}{Analysis and {Control} of {Flat} {Systems} by {Geometric} {Methods}}}, ser. Modellierung und {Regelung} komplexer dynamischer {Systeme}.\hskip 1em plus 0.5em minus 0.4em\relax Düren: Shaker Verlag, 2023, vol.~59.

\bibitem{pasillas-lepine_contact_2001}
W.~Pasillas-Lépine and W.~Respondek, ``Contact systems and corank one involutive subdistributions,'' \emph{Acta Applicandae Mathematica}, vol.~69, no.~2, pp. 105--128, Nov. 2001.

\bibitem{rudolph_examples_1998}
J.~Rudolph and E.~Delaleau, ``\BIBforeignlanguage{en}{Some examples and remarks on quasi-static feedback of generalized states},'' \emph{\BIBforeignlanguage{en}{Automatica}}, vol.~34, no.~8, pp. 993--999, Aug. 1998.

\end{thebibliography}

\appendices

\section{Sketch of the Proof of Theorem~\ref{thm:TFCCS+IC_0}}
\label{app:proof_thm_TFCCS_0}

    \emph{Necessity:} It is straightforward to verify that TF\(_0\) meets all conditions of Theorem~\ref{thm:TFCCS+IC_0}.
    
    \emph{Sufficiency:} To show sufficiency, we transform a system  \eqref{eq:ai_sys_m} that satisfies all conditions of Theorem~\ref{thm:TFCCS+IC_0} into TF$_0$ step by step using state and input transformations only.
    
    \emph{Step~1: Splitting off $\Sigma_\zeta$.} Consider a system of the form \eqref{eq:ai_sys_m} for which the conditions of Theorem~\ref{thm:TFCCS+IC_0} are met. With condition~1) of Theorem~\ref{thm:TFCCS+IC_0}, we obtain the sequence of involutive distributions \mbox{$\Dist_1 \crksub{m+1} \ldots \crksub{m+1} \Dist_{k^\zeta}$}. By straightening out these involutive distributions, we obtain
    \begin{equation}\label{eq:prf_thm_TFCCS_0_seq_1_straight}
        \begin{aligned}
        \Dist_1 &= \operatorname{span}\{\pad{\bar{\zeta}^{k^\zeta}}\}, \Dist_2 = \operatorname{span}\{\pad{\bar{\zeta}^{k^\zeta}}, \pad{\bar{\zeta}^{k^\zeta-1}}\}, \ldots\\
        D_{k^\zeta-1} &= \operatorname{span}\{ \pad{\bar{\zeta}^{k^\zeta}}, \ldots, \pad{\bar{\zeta}^{2}} \}, D_{k^\zeta} = \operatorname{span}\{ \pad{\bar \zeta} \} ,
        \end{aligned}
    \end{equation}
    where $\operatorname{span}\{\pad{\bar{\zeta}^{j}}\} = \operatorname{span}\{\pad{\bar{\zeta}^{j}_0}, \ldots, \pad{\bar{\zeta}^{j}_m}\}$ for any \mbox{$1 \le j \le k^\zeta$}. In these coordinates, \eqref{eq:ai_sys_m} decomposes into
    \begin{equation}\label{eq:prf_thm_TFCCS_0_sys_z}
    \begin{aligned}
        \Sigma_z: \hspace{9pt} & \zdot = f_z(z, \bar{\zeta}^1_0, \ldots, \bar{\zeta}^1_m)\; , \\
        \Sigma_\zeta: \hspace{9pt} & \hspace{0pt} \dot{\bar{\zeta}} = f_{\zeta}(z, \bar \zeta, u) ,
    \end{aligned}
    \end{equation}
    where $\Sigma_\zeta$ has a triangular structure and its top state variables \mbox{$\bar{\zeta}^1=(\bar{\zeta}^1_0, \ldots, \bar{\zeta}^1_m)$} act as inputs to $\Sigma_z$. Consequently, condition~2) of Theorem~\ref{thm:TFCCS+IC_0}, i.e., $\Dist_{k^\zeta} = \Cauchy{ \Dist_{k^\zeta+1} }$, implies 
    \begin{equation*}
        \Span{\pad{\bar{\zeta}^1}} = \Cauchy{\operatorname{span}\{\pad{\bar{\zeta}^1}, [\pad{\bar{\zeta}^1_0}, f_z], \ldots, [\pad{\bar{\zeta}^1_m}, f_z]\} }. 
    \end{equation*}
    Therefore, according to Lemma \ref{lem:control_affine}, there exists a state transformation \mbox{$\zeta^1 = (\zeta^1_0, \ldots, \zeta^1_m) = \Phi_{\zeta^1}(z,\bar{\zeta}^1)$} such that $\zeta^1$ appears as a control affine input in the upper subsystem of \eqref{eq:prf_thm_TFCCS_0_sys_z}. For notational convenience, we relabel $\zeta^1$ back to  $\bar \zeta^1$ and obtain
    \begin{equation*}\label{eq:prf_thm_0_sys_ai_z}
        \Sigma_z: \; \; \zdot = a_z(z) + b_{z,0}(z)\bar\zeta^1_0 + \cdots + b_{z,m}(z)\bar\zeta^1_m .
    \end{equation*} 
 
    \emph{Step~2: Decomposing $\Sigma_z$ into $\Sigma_\chi$ and $\Sigma_\xi$.} Next, we define the distribution $\E_z=\Span{b_{z,0},\ldots,b_{z,m}}$. By assumption, \mbox{$\E = \Dist_{k^\zeta} + \E_z$} fulfills condition~3) of Theorem~\ref{thm:TFCCS+IC_0}. Accordingly, $\E_z$ satisfies
    \begin{enumerate}[label=\alph*)]
        \item $\operatorname{rank}(\bar{\E_z}) = mk^\chi + 1$.
        \item $\operatorname{rank}(\E_z^{(k^\chi-2)}) = \operatorname{rank}(\bar{\E_z}) - m$ and there exists an involutive subdistribution $\Ld_z \crksub{1} \E_z^{(k^\chi-2)}$.
        \item The drift $a_z$ satisfies the compatibility conditions
        $[a_z, \Cauchy{\E_z^{(i)}}] \subset \E_z^{(i)}$, for all $1 \leq i \leq k^\chi-2$, and $[a_z, \Ld_z] \subset \bar{\E_z}$.
    \end{enumerate}
    From \cite{RespondekCanonicalContactSystems2002} it can be derived that the geometry of the subsystem described by $\E_z$ is given by the following diagram:
    \begin{equation*}
        \arraycolsep=0.5pt
        \begin{array}{ccccccccccc}
            \E_z & \underset{m}{\subset} & \cdots & \underset{m}{\subset} & \E_z^{(k^\chi-3)} & \underset{m}{\subset} & \E_z^{(k^\chi-2)} & \underset{m}{\subset} &  \E_z^{(k^\chi-1)} = \bar{\E}_z \\[-0.2em]
            \rotsubset &  & & & \rotsubset & & \rotsubset & & \\[0.9em]
            \Cauchy{\E_z^{(1)}} & \underset{m}{\subset} & \cdots & \underset{m}{\subset} & \Cauchy{\E_z^{(k^\chi-2)}} & \underset{m}{\subset} & \Ld_z .
        \end{array}
    \end{equation*}
    
    Item~4) of Theorem~\ref{thm:TFCCS+IC_0} implies that all distributions \mbox{$\F_{z,i+1} = \F_{z,i} + [a_z,\F_{z,i}]$}, \mbox{$\F_{z,0}=\bar \E_z$}, \mbox{$i=0,\ldots,k^\xi-1$}, are involutive. Consequently, we have the sequence of involutive distributions
    \begin{equation*}\label{eq:prf_thm_TFCCS_0_seq_3}
        \Cauchy{\E_z^{(1)}} \crksub{m} \cdots \crksub{m} \Cauchy{\E_z^{(k^\chi-2)}} \crksub{m} \Ld_z \hspace{-2pt} \crksub{m+1} \hspace{-2pt} \F_{z,0} \subset \cdots \subset \F_{z,k^\xi} .
    \end{equation*}
    Straightening out these involutive distributions yields
    \begin{equation}\label{eq:prf_thm_TFCCS_0_seq_step_2}
        \begin{aligned}
            \Cauchy{\E_z^{(1)}} & = \Span{\pad{\bar{\chi}^{k^\chi}}}, \\
            & \hspace{0.5em} \vdots \\
            \Cauchy{\E_z^{(k^\chi-2)}} & = \Span{\pad{\bar{\chi}^{k^\chi}}, \ldots, \pad{\bar{\chi}^{3}}}, \\
            \Ld_z & = \Span{\pad{\bar{\chi}^{k^\chi}}, \ldots, \pad{\bar{\chi}^{2}}}, \\
            \bar{\E_z} = \F_{z,0} & = \Span{\pad{\bar{\chi}}}, \ldots, {\F}_{z,k^\xi} = \Span{\pad{\bar{\chi}}, \pad{\bar{\xi}}} . \\
        \end{aligned}
    \end{equation}
    Note that $\operatorname{span}\{ \pad{\bar{\chi}}\} = \operatorname{span}\{ \pad{\bar{\chi}^{k^\chi}}, \ldots, \pad{\bar{\chi}^{1}} \} $, where $\operatorname{span}\{\pad{\bar{\chi}^j}\} = \operatorname{span}\{\pad{\bar{\chi}^j_1}, \ldots, \pad{\bar{\chi}^j_m}\}$ for $2 \le j \le k^\chi$, and $\operatorname{span}\{\pad{\bar{\chi}^1}\} = \operatorname{span}\{\pad{\bar{\chi}^1_0}, \ldots, \pad{\bar{\chi}^1_m}\}$. Using the adapted coordinates from \eqref{eq:prf_thm_TFCCS_0_seq_step_2} and given that $[a_z, \Ld_z] \subset \bar{\E}_z$, system \eqref{eq:prf_thm_TFCCS_0_sys_z} decomposes into the subsystems
    \begin{equation}\label{eq:prf_thm_TFCCS_0_system_step2}
        \begin{aligned}
            \Sigma_\xi \hspace{1pt}: & \hspace{4pt} \dot{\bar{\xi}}  =a_\xi( \bar{\xi}, \bar{\chi}^1_0,\ldots,\bar{\chi}^1_m ), \\
            \Sigma_\chi: & \; \dot{\bar{\chi}} =a_\chi(  \bar{\xi}, \bar{\chi}) + b_{\chi,0}(  \bar{\xi}, \bar{\chi})\bar\zeta^1_0 + \cdots + b_{\chi,m}(  \bar{\xi}, \bar{\chi})\bar\zeta^1_{m} , \\
            \Sigma_\zeta: & \hspace{4pt} \dot{\bar{\zeta}} = f_{\zeta}(\bar \xi, \bar \chi, \bar \zeta, u) .
        \end{aligned}
    \end{equation}
    
    \emph{Step~3: Transforming $\Sigma_\xi$ to Brunovský normal form.} Given the involutive sequence \mbox{$\F_{z,i+1}=\F_{z,i}+[a_z,\F_{z,i}]$}, \mbox{$i=0,\ldots,k^\xi-1$}, we know from the input-state linearizability problem, that we can transform $\Sigma_\xi$ into Brunovský normal form \eqref{eq:TFCCS+IC_xi} with inputs $(\chi^1_0,\ldots, \chi^1_m)$ using state transformations of the form $(\xi, (\chi^1_0, \ldots, \chi^1_m))=(\Phi_{\xi}(\bar{\xi}), \Phi_{\chi^1}(\bar{\xi}, \bar{\chi}^1_0,\ldots,\bar{\chi}^1_m)$, preserving \eqref{eq:prf_thm_TFCCS_0_seq_step_2} as well as the triangular structure of \eqref{eq:prf_thm_TFCCS_0_system_step2}.
    
    \emph{Step~4: Transforming $\Sigma_\chi$ to the contact form.}
    Due to the triangular structure, the top $\chi$-variables are no longer altered by the transformation and the $\xi$-variables play no role in the transformation of the subsystem $\Sigma_\chi$. The subsystem $\Sigma_\chi$ is control-affine and meets all conditions for conversion to the canonical contact form. Moreover, the input vector fields act only in the $\chi$-directions. For more details, we refer to the literature on contact forms with drift, in particular, \cite{LiMultiinputControlaffineSystems2016}. In the case $s=0$, i.e., for a control-affine system SFE to TF$_0$, the lowest equations in \eqref{eq:TFCCS+IC_eta} can actually be normalized to \mbox{$\dot{\chi}_j^{\,k^{\chi}}=\zeta_j^{1}$}, \mbox{$j=1,\ldots,m$}.

    \emph{Step~5: Transforming $\Sigma_\zeta$ to Brunovský normal form.} Straightening out the distributions \eqref{eq:prf_thm_TFCCS_0_seq_1_straight} in Step~1 puts the subsystem $\Sigma_\zeta$ in triangular form, which is preserved up to this point. Hence, similar to the standard static-feedback linearization problem, $\Sigma_\zeta$ can be transformed into \(m+1\) integrator chains, each of length $k^\zeta$.
    
\section{Sketch of the Proof of Theorem~\ref{thm:TFCCS+IC_1}}
\label{app:proof_thm_TFCCS_1}

\emph{Necessity:} Verifying that TF$_1$ satisfies all conditions of Theorem~\ref{thm:TFCCS+IC_1} is straightforward.

\emph{Sufficiency:} We transform a system \eqref{eq:ai_sys_m} satisfying Theorem~\ref{thm:TFCCS+IC_1} into TF\(_1\) using only state and static input transformations.

\emph{Step~1: Splitting off $\Sigma_\zeta$.} The goal of this step is to split off $\Sigma_\zeta$ such that $m+1$ state variables of $\Sigma_\zeta$ appear as control-affine inputs in the subsystem above denoted by $\Sigma_z$. By assumption, condition~1) holds. Unlike the TF$_0$ case, item~2) states $\Cauchy{\Dist_{k^\zeta+1}} \neq \Dist_{k^\zeta}$. Therefore, proceeding as in Step~1 of the sufficiency part of Appendix~\ref{app:proof_thm_TFCCS_0} does not yield a control-affine subsystem $\Sigma_z$. Consequently, we straighten out the truncated sequence $\Dist_1 \subset\cdots\subset \Dist_{k^\zeta-1}$ and obtain
\begin{equation*}
    \begin{aligned}
        \Dist_1 = \operatorname{span}\{\pad{\bar{\zeta}^{k^\zeta}}\}, \ldots, \Dist_{k^\zeta-1} &= \operatorname{span}\{\pad{\bar{\zeta}^{k^\zeta}}, \ldots, \pad{\bar{\zeta}^{2}}\}, 
    \end{aligned}
\end{equation*}
where \mbox{$\operatorname{span}\{\pad{\bar{\zeta}^j}\} = \operatorname{span}\{\pad{\bar{\zeta}^j_0}, \ldots, \pad{\bar{\zeta}^j_m}\}$} for any \mbox{$2 \le j \le k^\zeta$}. 

Arguing analogously to the sufficiency-proof of TF$_0$ and leveraging that $\operatorname{span}\{\pad{\bar{\zeta}^2}\} \subset \Cauchy{\Dist_{k^\zeta}}$, we obtain
\begin{equation*}
\begin{aligned}
    \zdot_1 & = a_{z_1}(z_1) + b_{z_1,0}(z_1)\bar{\zeta}^2_0 + \cdots + b_{z_1,m}(z_1)\bar{\zeta}^2_m,\\
    \dot{\bar{\zeta}}_2 & = f_{\bar \zeta_2}(z_1, \bar{\zeta}_2, u) \; ,
\end{aligned}
\end{equation*}
where $\bar \zeta_2=(\bar \zeta^{k^\zeta}, \ldots, \bar \zeta^{2})$ with $\bar \zeta^{j}=(\bar{\zeta}^{j}_{0},\ldots,\bar{\zeta}^{j}_{m})$ for any $2 \le j \le k^\zeta$. Given condition~1), the distribution \mbox{$\Dist_{z_1,1} = \operatorname{span}\{b_{z_1,0}, \ldots, b_{z_1,m}\}$} is involutive. Next, we define \mbox{$\Dist_{z_1,2} = \Dist_{z_1,1} + [a_{z_1}, \Dist_{z_1,1}]$}, which is by assumption non-involutive with $\Cauchy{\Dist_{z_1,2}} \neq \Dist_{z_1,1}$, see condition~2) of Theorem~\ref{thm:TFCCS+IC_1}.
Since condition~3) of Theorem~\ref{thm:TFCCS+IC_1} holds by assumption and given that \mbox{$\Dist_{k^\zeta} = \Dist_{z_1,1}+\Dist_{k^\zeta-1}$} and \mbox{$\Dist_{k^\zeta+1} = \Dist_{z_1,2}+\Dist_{k^\zeta-1}$}, it follows that there exist $m$ linearly independent vector fields 
\begin{equation*}
    c_i = \sum_{k=0}^{m} \alpha^k_ib_{z_1,k}, \quad i=1,\ldots,m,
\end{equation*}
such that \mbox{$\E_{z_1} = \Dist_{z_1,1} + \Span{[a_{z_1},c_1], \ldots, [a_{z_1},c_{m}]}$} satisfies \mbox{$\Cauchy{\E_{z_1}}=\Span{c_1,\ldots,c_{m}}$}. 

Then, apply a state transformation $\zeta^2=\Phi_{\zeta}(z_1,\bar \zeta^2)$ and relabel $\zeta^2$ back to $\bar \zeta^2$ such that 
\begin{equation}\label{eq:prf_thm_1_system_z1}
    \zdot_1 = a_{z_1}(z_1) + \bar b_{z_1,0}(z_1)\bar {\zeta}^2_0 + c_1(z_1)\bar\zeta^2_1 + \cdots + c_m(z_1)\bar\zeta^2_m .
\end{equation}
Straightening out the vector fields $c_1,\ldots,c_m$ yields \mbox{$\Dist_{z_1,1} = \operatorname{span}\{b_{z,0},\pad{\bar{\zeta}^1}\}$} with $\bar{\zeta}^1=(\bar \zeta^1_1, \ldots, \bar \zeta^1_m)$.
Consequently, we have split off $m$ states $(\bar{\zeta}^1_1,\ldots,\bar{\zeta}^1_m)$ from \eqref{eq:prf_thm_1_system_z1}. Without loss of generality, we assume that $b_{z,0}$ is associated to the input $\bar\zeta^2_0$. After these transformations, the complete system is of the form
\begin{equation}\label{eq:prf_thm_TFCCS_1_sys_z}
\begin{aligned}
    \Sigma_z: \hspace{9pt} & \zdot = a_z(z,\bar \zeta^1) + b_{z,0}(z,\bar \zeta^1)\bar \zeta^2_0 = f_z(z, \bar \zeta^1, \bar \zeta^2_0), \\
    \Sigma_\zeta: \hspace{9pt} & \hspace{0pt} \dot{\bar{\zeta}} = f_{\zeta}(z, \bar \zeta, u) \; ,
\end{aligned}
\end{equation}
with $\Sigma_\zeta$ in triangular form. Its top state variables \mbox{$(\bar{\zeta}^1, \bar{\zeta}^2_0)$} act as inputs to $\Sigma_z$. Based on $\operatorname{span}\{ \pad{\bar \zeta ^1} \} = \Cauchy{\E_{z_1}} = \Cauchy{\operatorname{span}\{\pad{\bar \zeta^1}, b_{z,0}, [\pad{\bar\zeta^1}, a_z] \}}$ it can be shown that $\Sigma_z$ can be completely brought into an affine input form. In particular, this is possible by a transformation of the form \mbox{$(\hat\zeta^1,\hat\zeta_0^2)=(\Phi_{\hat\zeta^1}(z,\bar\zeta^1),\Phi_{\hat\zeta^2_0}(z,\bar\zeta^1,\bar\zeta_0^2))$}, which by construction preserves the triangular structure of $\Sigma_\zeta$ established so far.
Consequently, $\Sigma_z$ allows the control-affine representation
\begin{equation}\label{eq:prf_thm_1_sys_ai_z}
    \zdot = a_z(z) + b_{z,0}(z)\bar\zeta^2_0 + b_{z,1}(z)\bar{\zeta}^1_1 + \cdots +  b_{z,m}(z)\bar{\zeta}^1_m \; ,
\end{equation}
where, for consistency, we returned to the $(z,\bar \zeta)$-notation used in \eqref{eq:prf_thm_TFCCS_1_sys_z}. By defining the input distribution \mbox{$\E_z =\Span{ b_{z,0},\ldots, b_{z,m}}$}, we return to the starting point from Step~2 of the sufficiency part in Appendix \ref{app:proof_thm_TFCCS_0}. Thus, for the remaining subsystem \eqref{eq:prf_thm_1_sys_ai_z}, the TF$_1$ case ($s=1$) reduces to the TF$_0$ case ($s=0$). The rest of the proof is identical to Steps 2–4 in the proof of Theorem~\ref{thm:TFCCS+IC_0}, with a minor difference at the end of Step~4. Here, the lowest state equations corresponding to $\chi^{k^\chi}_j,j=1,\ldots,m$, cannot generally be normalized to $\dot{\chi}^{k^\chi}_j=\zeta^1_j,  j=1,\ldots,m$, since $\zeta^1_j$ can have a relative degree of $k^\zeta-1$ which in turn would not comply with the triangular form \eqref{eq:TFCCS+IC_q}. 

\begin{rem}[Case \(s=1\) vs.\ \(s\ge 2\)]

\emph{Case \(s=1\).} It suffices to identify the states belonging to \(\Sigma_\zeta\) and to verify that the remaining subsystem is locally static-feedback equivalent to the (canonical) contact form with appended integrator chains. The reason is that, in adapted coordinates, the state \(\chi^0\) associated with the distinguished input \(w^0\) and the states \((\zeta^1_1,\ldots,\zeta^1_m)\) acting as inputs of \(\Sigma_\chi\) share the same relative degree \(k^\zeta\). Consequently, any state transformation \mbox{$(\tilde{\zeta}^1_1,\ldots,\tilde{\zeta}^1_m)=\Phi(\chi,\zeta^1_1,\ldots,\zeta^1_m)$} is compatible with the triangular structure \eqref{eq:TFCCS+IC_q}.

\emph{Case \(s\ge 2\).} After identifying \(\Sigma_\zeta\), one can still analyze the remaining subsystem separately and test static-feedback equivalence to the contact form with integrator chains. If the test succeeds, the overall system can be rendered into a structurally flat triangular form similar to \eqref{eq:TFCCS+IC_q}, with the difference that the lower subsystem $\Sigma_\zeta$ need not be a pure stack of integrator chains. In general, the derivatives \(\dot{\zeta}^{\,i}_j\), for any \(i\in\{k^\zeta-s,\ldots,k^\zeta-1\}\) and  \(j\in\{1,\ldots,m\}\), may depend explicitly on the distinguished input \(w^0\).
\end{rem}

\newpage

\end{document}